\RequirePackage{ifpdf}
\ifpdf 
\documentclass[pdftex]{sigma}
\else
\documentclass{sigma}
\fi

\usepackage{mathrsfs}

\begin{document}

\allowdisplaybreaks

\renewcommand{\thefootnote}{$\star$}

\renewcommand{\PaperNumber}{038}

\FirstPageHeading

\ShortArticleName{Elliptic Hypergeometric Solutions to Elliptic Dif\/ference Equations}

\ArticleName{Elliptic Hypergeometric Solutions\\ to Elliptic Dif\/ference Equations\footnote{This paper is a contribution to the Proceedings of the Workshop ``Elliptic Integrable Systems, Isomonodromy Problems, and Hypergeometric Functions'' (July 21--25, 2008, MPIM, Bonn, Germany). The full collection
is available at
\href{http://www.emis.de/journals/SIGMA/Elliptic-Integrable-Systems.html}{http://www.emis.de/journals/SIGMA/Elliptic-Integrable-Systems.html}}}

\Author{Alphonse P. MAGNUS}

\AuthorNameForHeading{A.P. Magnus}

\Address{Universit\'e catholique   de Louvain,
           Institut math\'ematique,\\ 2 Chemin du Cyclotron,
            B-1348  Louvain-La-Neuve, Belgium}

\Email{\href{mailto:alphonse.magnus@uclouvain.be}{alphonse.magnus@uclouvain.be}}
\URLaddress{\url{http://perso.uclouvain.be/alphonse.magnus/}}

\ArticleDates{Received December 01, 2008, in f\/inal form March 20,
2009; Published online March 27, 2009}

\Abstract{It is shown how to def\/ine dif\/ference
equations on particular lattices $\{x_n\}$, $n\in\mathbb{Z}$, made of
 values of an elliptic function at a sequence of arguments in
arithmetic progression (\emph{elliptic lattice}). Solutions to special
dif\/ference equations have remarkable
simple interpolatory expansions.
Only linear dif\/ference equations of f\/irst order
are considered here.}

\Keywords{elliptic dif\/ference equations;  elliptic hypergeometric expansions}

\Classification{39A70; 41A20}

  \begin{flushright}  \sl Nacht und St\"urme werden Licht\\ Choral Fantasy, Op. 80
      \end{flushright}

\section[Difference equations on elliptic lattices]{Dif\/ference equations on elliptic lattices}

\subsection[The difference operator]{The dif\/ference operator}

We consider functional equations involving the dif\/ference operator
\begin{gather}\label{diffop}
 (\mathcal{D} f)(x) = \frac{ f(\psi(x))- f(\varphi(x)) }{\psi(x)-\varphi(x) }.
\end{gather}

Most
instances \cite{MilneT} are $(\varphi(x),\psi(x))= (x,x+h)$, or the more
symmetric $(x-h/2,x+h/2)$, or also $(x,qx)$ in $q$-dif\/ference
equations \cite{GrH,Koe,Ko}. Recently, more complicated forms  $(r(x)-\sqrt{s(x)},
r(x)+\sqrt{s(x)})$ have appeared \cite{AAR,AW,Koe,Ko,Segovia,magsnul,NS,Niki,MagLuminy2007}, where
$r$ and $s$ are rational functions.

This latter trend will be examined here: we need, for each $x$,
two values $f(\varphi(x))$ and $f(\psi(x))$ for $f$.
A f\/irst-order dif\/ference equation is
\[
\mathcal{F}\left(x,f(\varphi(x)),f(\psi(x))\right)=0, \qquad {\rm or} \qquad
f(\varphi(x))-f(\psi(x))=
\mathcal{G}\left(x,f(\varphi(x)),f(\psi(x))\right)\] if we want to
emphasize the dif\/ference of $f$. There is of course some freedom
in this latter writing. Only \emph{symmetric} forms in $\varphi$
and $\psi$ will be considered here:
\begin{gather*}
 (\mathcal{D} f)(x) = \mathscr{F}\left(x,f(\varphi(x)),f(\psi(x))\right),
\end{gather*}
where $\mathcal{D}$ is the divided dif\/ference operator \eqref{diffop}
 and where $\mathscr{F}$ is a symmetric function of its two last arguments.

For instance, a linear dif\/ference equation of  f\/irst order may be written as
\[   a(x) f(\varphi(x)) + b(x) f(\psi(x)) +c(x)=0,
\]
as  well  as
\[  \alpha(x) (\mathcal{D} f)(x) = \beta(x)
[f(\varphi(x)) + f(\psi(x))] +\gamma(x),\]
with
$\alpha(x)=[b(x)-a(x)][\psi(x)-\varphi(x)]/2$,
$\beta(x)=-[a(x)+b(x)]/2$, and $\gamma(x)=-c(x)$.

  The simplest choice
 for $\varphi$ and $\psi$ is to take the two determinations of an
algebraic function of degree 2, i.e., the two $y$-roots of
\begin{subequations}
\begin{gather}\label{Fy}
    F(x,y)=X_0(x)+X_1(x) y +X_2(x) y^2 =0,
\end{gather}
 where $X_0$, $X_1$, and $X_2$ are rational functions.

Note that the sum and the product of $\varphi$ and $\psi$ are
the rational functions
\begin{gather}\label{sumprod}
 \varphi+\psi=-X_1/X_2,\qquad \varphi\psi=X_0/X_2.
\end{gather}

\subsection{\label{latti}The corresponding lattice, or grid}

 Dif\/ference equations must allow the recovery of $f$ on a whole set
 of points. An initial-value problem for a f\/irst order dif\/ference equation
 starts with a value for $f(y_0)$ at $x=x_0$, where $y_0$ is one root
 of \eqref{Fy} at $x=x_0$. The dif\/ference equation at $x=x_0$ relates then
 $f(y_0)$ to $f(y_1)$, where $y_1$ is the second root of \eqref{Fy} at $x_0$.
 We need $x_1$ such that $y_1$ is one of the two roots of~\eqref{Fy} at~$x_1$,
so for one of the roots of $F(x,y_1)=0$ which is not~$x_0$. Here again, the
simplest case is when $F$ is of degree~2 in~$x$:
 \begin{gather}\label{Fx}
   F(x,y)=Y_0(y)+Y_1(y) x +Y_2(y) x^2 =0.
\end{gather}
\end{subequations}
Both forms \eqref{Fy} and \eqref{Fx} hold simultaneously
when $F$ is \emph{\textbf{biquadratic}}:
\begin{gather}\label{Fbiq}
 F(x,y)= \sum_{i=0}^2 \sum_{j=0}^2 c_{i,j} x^i y^j.
\end{gather}

The construction where successive points on the curve $F(x,y)=0$
are $(x_n, y_n)$, $(x_n,y_{n+1})$, $(x_{n+1},y_{n+1})$, is called
``T-algorithm'' in \cite[Theorem~6]{SpZ2007}, see also the Fritz
John's algorithm in~\cite{BursZhed1,BursZhed,BursZhed09}. The sequence
$\{x_n\}$ is then an instance of \emph{\textbf{elliptic}} lattice,
or grid.

Of course, the sequence $\{y_n\}$ is elliptic too, $x_n$ and $y_n$
have elliptic functions representations
\begin{gather} \label{ell12}
x_n=\mathcal{E}_1(t_0+nh), \qquad  y_n=\mathcal{E}_2(t_0+nh),
\end{gather}
where $(x=\mathcal{E}_1(t), y=\mathcal{E}_2(t))$ is a parametric
representation of the biquadratic curve $F(x,y)=0$ with the $F$ of~\eqref{Fbiq}.

Note that the names of the $x$- and $y$-lattices are sometimes inverted, as
in \cite[equation~(1.2)]{SpZ2007}

As $y_n$ and $y_{n+1}$ are the two roots in $t$ of $F(x_n,t)= X_0(x_n)+X_1(x_n)t
+X_2(x_n)t^2=0$, useful identities are
\begin{gather*} 
y_n+y_{n+1} =-\frac{X_1(x_n)}{X_2(x_n)},\qquad  y_n y_{n+1} =\frac{X_0(x_n)}{X_2(x_n)},
\end{gather*}
from \eqref{sumprod}, and the direct formula
\begin{gather*} 
y_n \  \text{and}  \ y_{n+1} = \frac{ -X_1(x_n) \pm \sqrt{P(x_n)}  }{2 X_2(x_n) },
\end{gather*}
where
\begin{gather*}
P=X_1^2-4X_0 X_2\end{gather*} is a polynomial of degree 4.

Also, as $x_{n+1}$ and $x_n$ are the two roots in $t$ of $F(t,y_{n+1})=0$,
\begin{gather*} 
x_n+x_{n+1} =-\frac{Y_1(y_{n+1})}{Y_2(y_{n+1})},\qquad  x_n x_{n+1} =\frac{Y_0(y_{n+1})}{Y_2(y_{n+1})}.
\end{gather*}

As the operators considered here are symmetric in $\varphi(x)$ and
$\psi(x)$, we do not need to def\/ine precisely what $\varphi$ and
$\psi$ are, i.e., we only need to know the pair $(\varphi,\psi)$,
and not the ordered pair. However, once a starting point
$(x_0,y_0)$ is chosen, it will be convenient to def\/ine
$\varphi(x_n)=y_n$ and $\psi(x_n)=y_{n+1}$, $n\in\mathbb{Z}$.

\emph{Special cases.} We already encountered the usual dif\/ference operators
$(\varphi(x),\psi(x))=(x,x+h)$ or $(x-h,x)$ or $(x-h/2,x+h/2)$ corresponding to
$X_2(x)\equiv 1$, $X_1$ of degree 1, $X_0$ of degree 2 with $P=X_1^2-4X_0X_2$ of degree 0.
For the geometric dif\/ference operator, $P$ is the square of a~f\/irst degree polynomial.
For the Askey--Wilson operator \cite{AAR,AW,Ism2005,Koe,Segovia,magsnul}, $P$ is an arbitrary second degree
polynomial.

The formulas for the sequences $x_n$ and $y_n$  are in these three special cases
\begin{gather*}
(x_n,y_n)= (x_0+nh, y_0+nh);\qquad (a+bq^n,u+v q^n);\nonumber\\ (a+bq^n+cq^{-n},u+vq^n+wq^{-n}). 
\end{gather*}

\subsection[Difference of a rational function]{Dif\/ference of a rational function}

From \eqref{sumprod},
when the divided dif\/ference operator $\mathcal{D}$ of \eqref{diffop} is applied to
a rational function, the result is still a rational function.

The dif\/ference operator applied to a simple rational function is of
special interest.

Let $f(x)=\frac{1}{x-A}$, then
\begin{gather*}
 \mathcal{D} \frac{1}{x-A}  = \frac{1}{\psi(x)-\varphi(x)}\left[
   \frac{1}{\psi(x)-A} -\frac{1}{\varphi(x)-A} \right] =
 -\frac{1}{(\psi(x)-A)(\varphi(x)-A)} \\
\phantom{\mathcal{D} \frac{1}{x-A}}{} = -\frac{X_2(x)}{X_0(x)+A X_1(x)+A^2 X_2(x)},
\end{gather*}
and let
$\{(x'_n,y'_n), (x'_n,y'_{n+1})\}$ be the  elliptic sequence on the
biquadratic curve $F(x,y)=0$ such that $y'_0=A$, then
\begin{gather} \label{simplerat}
\mathcal{D} \frac{1}{x-A}
=-\frac{X_2(x)}{Y_2(A)(x-x'_0)(x-x'_{-1})},
\end{gather}
as the denominator is $F(x,A)$, and the two $x$-roots of
$F(x,A)=F(x,y'_0)=0$ are $x'_0$ and $x'_{-1}$, from the opening
discussion of Section~\ref{latti}.

The $\mathcal{D}$ operator applied to a general rational function
yields a rational function with the factor $X_2$. It seems
sometimes f\/it to def\/ine a dif\/ference operator as our $\mathcal{D}$
divided by $X_2$, as by V.P.~Spiridonov and A.S.~Zhedanov in
Section~6 of~\cite{SpZ1}. See also Section~2 of \cite{SpZ2007}.

A general rational function is generically a sum of simple
rational functions of type \eqref{simplerat}, say, $1/(x-A)$,
$1/(x-B)$, etc. The dif\/ference has poles at $x'_0$ and $x'_{-1}$,
also at $x''_0$ and $x''_{-1}$ if $B=y''_0$, etc., so that the degree
of $\mathcal{D}f$ is usually twice the degree of $f$.
However, the dif\/ference of a rational function of denominator
$(x-y'_0)(x-y'_1)\cdots (x-y'_n)$, $\mathcal{D}f$~has no other
poles than $x'_{-1}, x'_0,\dots, x'_n$.
This is also discussed in \cite{SpZ1,SpZ2007}.

So, let $\{(x_n,y_n), (x_n,y_{n+1})\}$ be a f\/irst elliptic sequence on
the biquadratic curve $F(x,y)=0$, and $\{(x'_n,y'_n),
(x'_n,y'_{n+1})\}$ be another  elliptic sequence on the same
curve. The two sequences have the same formula \eqref{ell12}, but
with dif\/ferent starting values $t_0$ and $t'_0$.

Now, let  \[ \mathcal{X}_n(x)=  \frac{(x-x_0)\cdots (x-x_{n-1}) }{ (x-x'_1)\cdots (x-x'_n)}\qquad {\rm and}
\qquad \mathcal{Y}_n(x)=\frac{(x-y_0)\cdots (x-y_{n-1}) }{(x-y'_1)\cdots (x-y'_n)}.
\]
See that \begin{gather*} 
 \mathcal{D}\mathcal{Y}_n(x) =
      C_n X_2(x) \frac{\mathcal{X}_{n-1}(x)}{(x-x'_0)(x-x'_n)}.
\end{gather*}

Indeed, $(\varphi(x)-y_0)(\varphi(x)-y_1)\cdots  (\varphi(x)-y_{n-1})$ and
  $(\psi(x)-y_0)(\psi(x)-y_1)\cdots  (\psi(x)-y_{n-1})$ both vanish at
$x=x_0,x_1,\dots, x_{n-2}$;
$(\varphi(x)-y'_1)(\varphi(x)-y'_2)\cdots  (\varphi(x)-y'_n)$ vanishes at
$x=x'_1,\dots, x'_n$, whereas
$(\psi(x)-y'_1)(\psi(x)-y'_2)\cdots  (\psi(x)-y'_n)$ vanishes at
$x=x'_0,\dots, x'_{n-1}$.

Simple fractions give
\[
\mathcal{D}\frac{1}{x-y'_k}=
 -\frac{X_2(x)}{Y_2(y'_k)(x-x'_{k-1})(x-x'_k)},
 \] as seen earlier in \eqref{simplerat}.

The constant $C_n$ is found through particular values of $x$, either $x_{-1}$,
where $\mathcal{Y}_n(\psi(x))=0$ but $\mathcal{Y}_n(\varphi(x))\neq 0$, or $x_{n-1}$,
where $\mathcal{Y}_n(\varphi(x))=0$ but $\mathcal{Y}_n(\psi(x))\neq 0$:
\begin{subequations}\label{Cx}
\begin{gather}
C_n = -\frac{\mathcal{Y}_n(\varphi(x_{-1})=y_{-1})
                (x_{-1}-x'_0)(x_{-1}-x'_n)   }
                  {(y_0-y_{-1})X_2(x_{-1})\mathcal{X}_{n-1}(x_{-1})},\label{Cxm1}
\\
   C_n = \frac{\mathcal{Y}_n(\psi(x_{n-1})=y_n)
                    (x_{n-1}-x'_0)(x_{n-1}-x'_n)  }
         {(y_n-y_{n-1})X_2(x_{n-1})  \mathcal{X}_{n-1}(x_{n-1})}\label{Cxn1}
\end{gather}
(of course, $C_0=0$).
Or through residues at $x'_0$, where $\mathcal{Y}_n(\psi(x))=\infty$, or $x'_n$
where $\mathcal{Y}_n(\varphi(x))=\infty$,
\begin{gather}\label{Cxp0}
C_n = \frac{(y'_1-y_0)\cdots (y'_1-y_{n-1})}
             {\frac{d\psi(x'_0)}{dx} (y'_1-y'_2)\cdots (y'_1-y'_n) }
              \frac{x'_0-x'_n}{(y'_1-y'_0)X_2(x'_0)\mathcal{X}_{n-1}(x'_0)},
\\
\label{Cxpn}
C_n=-
\frac{(y'_n-y_0)\cdots (y'_n-y_{n-1})}
             {(y'_n-y'_1)\cdots (y'_n-y'_{n-1})\frac{d\varphi(x'_n)}{dx} }
             \frac{x'_n-x'_0}{(y'_{n+1}-y'_n)X_2(x'_n)\mathcal{X}_{n-1}(x'_n)}.
\end{gather}
\end{subequations}

We shall also need the operator $\mathcal{M}$ def\/ined as
\begin{gather*} 
(\mathcal{M} f)(x)=[f(\varphi(x))+f(\psi(x))]/2,
\end{gather*}
which sends rational functions to rational functions too, usually of
double degree, but without particular factor.

With this operator
 $\mathcal{M}$,
\begin{gather*}
 2(\mathcal{M}  \mathcal{Y}_n)(x)=
 \frac{(\varphi(x)-y_0)(\varphi(x)-y_1)\cdots  (\varphi(x)-y_{n-1}) }
       {(\varphi(x)-y'_1)(\varphi(x)-y'_2)\cdots  (\varphi(x)-y'_n)}\\
\phantom{2(\mathcal{M}  \mathcal{Y}_n)(x)=}{}
   +\frac{ (\psi(x)-y_0)(\psi(x)-y_1)\cdots  (\psi(x)-y_{n-1})}{
             (\psi(x)-y'_1)(\psi(x)-y'_2)\cdots  (\psi(x)-y'_n)} \\
\phantom{2(\mathcal{M}  \mathcal{Y}_n)(x)}{}  =  2D_n(x) \frac{(x-x_0)(x-x_1)\cdots  (x-x_{n-2}) }
                {(x-x'_0)(x-x'_1)\cdots  (x-x'_n) }
 = 2 D_n(x)\frac{\mathcal{X}_{n-1}(x)}{(x-x'_0)(x-x'_n)} ,
\end{gather*}
where $D_n$ is a polynomial of degree 2.

Interesting values are found at the same point as in \eqref{Cx}:
\begin{subequations}\label{Dx}
\begin{gather}\label{Dxm1}
D_n(x_{-1}) = -\frac{C_n X_2(x_{-1}) (y_0-y_{-1})}{2},
\\
\label{Dxn1}
D_n(x_{n-1}) = \frac{C_n X_2(x_{n-1}) (y_n-y_{n-1})}{2},
\\
\label{Dxp0}
D_n(x'_0) = \frac{C_n X_2(x'_0) (y'_1-y'_0)}{2},
\\
\label{Dxpn}
D_n(x'_n) = -\frac{C_n X_2(x'_n) (y'_{n+1}-y'_n)}{2},
\end{gather}
\end{subequations}
when $n>0$. Of course, $D_0=1$.

\section{Elliptic hypergeometric expansions}

Let us consider expansions of the form
\begin{gather*} 
\sum_{k=0}^\infty   \prod_j (z_0^{(j)})^{\pm 1}
(z_1^{(j)})^{\pm 1} \cdots (z_k^{(j)})^{\pm 1},
\end{gather*}
where $z_k^{(j)}$ is a combination $a_j x_k^{(j)} +b_j$ or $a_j y_k^{(j)}+b_j$,
$\{\dots (x_k^{(j)},y_k^{(j)}), (x_k^{(j)},y_{k+1}^{(j)}),\dots\}$ being
elliptic lattices, or grids, related to a biquadratic curve \eqref{Fbiq}, the same
curve for each $j$.

We certainly recover at least a special case of current elliptic hypergeometric
expansions, as introduced in \cite{BursZhed1,BursZhed,Spir2008,SpZ1,SpZ2007}.

\subsection{Rational interpolatory elliptic expansions}

Rational interpolants of some function $f$ at $y_0$, $y_1,\dots,$
with poles at $y'_1, y'_2,\dots,$ are successive sums
\begin{gather}
c_0=f(y_0), \qquad c_0+c_1 \frac{x-y_0}{x-y'_1}, \qquad c_0+c_1 \frac{x-y_0}{x-y'_1}
  +c_2 \frac{(x-y_0)(x-y_1)}{(x-y'_1)(x-y'_2)},  \qquad  \dots, \nonumber\\ \sum_{k=0}^\infty c_k \mathcal{Y}_k(x).\label{expan}
\end{gather}

If, by chance, $c_k$ shows a similar form of ratio of products, we see special
cases of hypergeometric expansions! This will happen when one expands solutions
of dif\/ference equations which are simple enough. Putting the expansion
in the dif\/ference equation results in recurrence relations for~$c_k$, and we look for cases when this recurrence relation only involves
two terms~$c_k$ and~$c_{k+1}$.

\subsection[Linear 1${}^{st}$ order difference equations]{Linear 1${}^{\text{st}}$ order dif\/ference equations}

\begin{gather} \label{difflin}
a(x)(\mathcal{D} f)(x) = c(x) (\mathcal{M} f)(x)+d(x)
\end{gather}
Where is $b$? The full f\/lexibility of f\/irst order dif\/ference equations is
achieved with the
Riccati form \cite{MagLuminy2007}
\[
a(x)(\mathcal{D} f)(x) = b(x) f(\varphi(x))f(\psi(x))+c(x) [f(\varphi(x))+f(\psi(x))]+d(x)
\]
but only linear equations will be considered here.
However, \eqref{difflin} already allows elliptic exponentials ($c(x)\equiv a(x)$) or
logarithms ($c(x)\equiv 0$).

We now try to expand a solution to \eqref{difflin} as an interpolatory
series. If the initial condition is~$f(y_0)$ at~$x=x_0$, the dif\/ference equation allows
to f\/ind
\[
f(y_1)=  \frac{[a(x_0)/(y_1-y_0)+c(x_0)/2]f(y_0)+d(x_0)}{a(x_0)/(y_1-y_0)-c(x_0)/2} ,
\qquad f(y_2), \qquad \dots.
 \] This works f\/ine if no division by zero
is encountered.
Let us call $x'_0$ one of the roots of the algebraic
equation
\begin{gather} \label{xp0}
  \frac{a(x)}{\psi(x)-\varphi(x)}-\frac{c(x)}{2}=0,  \qquad  \text{at} \quad x=x'_0
\end{gather}
 and let, as usual,
$\psi(x'_0)=y'_1$, $\varphi(x'_0)=y'_0$. This shows that $y'_1$ is
a singular point of $f$, as trying to compute $f(y'_1)$ from
$f(y'_0)$ requires a division by zero. Then $y'_2,y'_3,\dots$
are poles as well.   That's why the expansion in \eqref{expan}
starts with poles at $y'_1,y'_2,\dots$ We also see that such
expansions represent meromorphic functions with a natural boundary
made of poles. At least, if the poles are spread on a curve, this
will be discussed in Section~\ref{conv}.

We also manage to have the initial value $f(y_0)$ completely
determined by the equation, i.e., independent of $f(y_{-1})$, so,
considering
\[
f(y_0)=
\frac{[a(x_{-1})/(y_0-y_{-1})+c(x_{-1})/2]f(y_{-1})+d(x_{-1})}{a(x_{-1})/(y_0-y_{-1})-c(x_{-1})/2},
\]
we ask $x_{-1}$ to be a root of
\begin{gather} \label{xm1}
  \frac{a(x)}{\psi(x)-\varphi(x)}+\frac{c(x)}{2}=0,  \qquad   \text{at} \quad  x=x_{-1}.
\end{gather}

Finally, we shall need the polynomials $c$ and $d$ to be of degree 3,
with $X_2$ as factor:
\begin{gather} \label{factX2}
c(x)=(\beta x+\gamma)X_2(x),\qquad  d(x)=(\delta x+\epsilon)X_2(x).
\end{gather}

We now have enough information for understanding the

\begin{theorem}\label{theorem1}
The difference equation \eqref{difflin} on the elliptic
lattice $F(x_n,y_n)=0$ of \eqref{Fy}--\eqref{Fbiq}, where $a$, $c$, and $d$ are polynomials
of degree $\leqslant 3$, $X_2$ being a factor of $c$ and $d$ as in \eqref{factX2}, has a solution with the formal
expansion \eqref{expan}, where $x_{-1}$ is a root of \eqref{xm1} and $x'_0$
is a root of \eqref{xp0}, with
\begin{gather*}
c_0=f(y_0)= \frac{d(x_{-1})}{a(x_{-1})/(y_0-y_{-1})-c(x_{-1})/2}=-\frac{d(x_{-1})}{c(x_{-1})}
=-\frac{\delta x_{-1}+\epsilon}{\beta x_{-1}+\gamma},\\
c_1= \frac{(\delta+\beta c_0)(x_0-x'_1)}{C_1(a(x_0)-c(x_0)(y_1-y_0)/2)}=
      \frac{(\gamma\delta-\beta\epsilon)(y_1-y'_1)X_2(x'_0)}{(y_1-y'_0)(x_0-x'_0)[a(x_0)-c(x_0)(y_1-y_0)/2]},
\end{gather*}
 and when $n\geqslant 1$,
\begin{gather}
c_n = c_1 \frac{C_1}{ x'_1-x_0 }  \frac{ x'_n-x_{n-1} }{C_n} \prod_{k=1}^{n-1}
     \frac{a(x'_k) +c(x'_k) (y'_{k+1}-y'_k)/2}{ a(x_k) -c(x_k) (y_{k+1}-y_k)/2}
  \frac{(x_k-x_{-1})(x_k-x'_0)}{(x'_k-x_{-1})(x'_k-x'_0)}\nonumber\\
\phantom{c_n}{} = -c_1 \frac{C_1}{ x'_1\!-x_0 }  (x'_n\!-x_{n-1})
     \frac{  (y_{-1}\!-y'_1)\cdots (y_{-1}\!-y'_{n-1}) X_2(x_{-1}) (x_{-1}\!-x_0)\cdots(x_{-1}\!-x_{n-2})  }
           {  (y_{-1}\!-y_1) \cdots (y_{-1}\!-y_{n-2})  (x_{-1}\!-x'_0)\cdots (x_{-1}\!-x'_n)  } \nonumber\\
\phantom{c_n=}{} \times \prod_{k=0}^{n-1}\frac{a(x'_k) +c(x'_k) (y'_{k+1}-y'_k)/2}{ a(x_k) -c(x_k) (y_{k+1}-y_k)/2}\;%
              \frac{(x_k-x_{-1})(x_k-x'_0)}{(x'_k-x_{-1})(x'_k-x'_0)}. \label{cn}
\end{gather}
\end{theorem}

\begin{proof}
Put the expansion \eqref{expan} in
\begin{gather*}
d(x)=a(x)\mathcal{D}f(x)-c(x)\mathcal{M}f(x)  =
\sum_0^\infty c_n\left[a(x)\mathcal{D}\mathcal{Y}_n(x)-c(x) (\mathcal{M} \mathcal{Y}_n(x)\right]
\\
\phantom{d(x)}{}  =-c_0  c(x)
+\sum_1^\infty c_n\left[ a(x)C_n X_2(x)
   -c(x) D_n(x) \right]\frac{\mathcal{X}_{n-1}(x)}{(x-x'_0)(x-x'_n)} .
\end{gather*}
The polynomial  $a(x)C_n X_2(x) -c(x) D_n(x) = [a(x)C_n -(\beta x+\gamma)D_n(x)]X_2(x)$ already has $X_2$ as factor
from \eqref{factX2}. A factor of degree
$\leqslant 3$ remains. Complete factoring follows:

at $x_{-1}$, from \eqref{Dxm1} and \eqref{xm1},
\[
a(x)C_n X_2(x) -c(x) D_n(x) = C_n X_2(x_{-1})[a(x_{-1}) +(y_0-y_{-1})c(x_{-1})/2]=0;
\]

at $x'_0$, from \eqref{Dxp0} and \eqref{xp0},
\[
a(x)C_nX_2(x)-c(x)D_n(x) =C_n X_2(x'_0)[a(x'_0)-(y'_1-y'_0)c(x'_0)/2]=0.
\]
Therefore we have three factors of f\/irst degree
\[
a(x)C_n X_2(x) -c(x) D_n(x)= X_2(x) (x-x_{-1})(x-x'_0)
[\xi_n(x-x_{n-1})+\eta_n(x-x'_n)],
\]
where from \eqref{Dxpn}
\begin{gather*} 
\xi_n = \frac{ a(x'_n)C_n X_2(x'_n) -c(x'_n) D_n(x'_n)}{ X_2(x'_n) (x'_n-x_{-1})(x'_n-x'_0)(x'_n-x_{n-1}) }
 = C_n \frac{ a(x'_n) +c(x'_n) (y'_{n+1}-y'_n)/2}{ (x'_n-x_{-1})(x'_n-x'_0)(x'_n-x_{n-1}) },
\end{gather*} and from \eqref{Dxn1}
\begin{gather*}
\eta_n = \frac{ a(x_{n-1})C_n X_2(x_{n-1}) -c(x_{n-1}) D_n(x_{n-1})}{ X_2(x_{n-1}) (x_{n-1}-x_{-1})(x_{n-1}-x'_0)(x_{n-1}-x'_n) }\nonumber\\
\phantom{\eta_n}{}
=C_n \frac{ a(x_{n-1}) -c(x_{n-1}) (y_n-y_{n-1})/2}{ (x_{n-1}-x_{-1})(x_{n-1}-x'_0)(x_{n-1}-x'_n) }.
\end{gather*}
Next,
\begin{gather*}
0 =a(x)\mathcal{D}f(x)-c(x)\mathcal{M}f(x) -d(x)\\
\phantom{0}{}  =-c_0 c(x)-d(x)
 +\sum_1^\infty c_n\left[ a(x)C_n X_2(x)
   -c(x) D_n(x) \right]\frac{\mathcal{X}_{n-1}(x)}{(x-x'_0)(x-x'_n)} \\
\phantom{0}{} =-c_0 c(x)-d(x)\\
\phantom{0=}{}
+\sum_1^\infty c_n X_2(x)\left[\xi_n(x-x_{n-1})+\eta_n(x-x'_n)
 \right]\frac{(x-x_{-1})(x-x_0)\cdots (x-x_{n-2})}{(x-x'_1)\cdots(x-x'_n)} \\
\phantom{0}{} =-c_0 c(x)-d(x)
+X_2(x)\sum_1^\infty c_n \xi_n\frac{(x-x_{-1})(x-x_0)\cdots (x-x_{n-2})(x-x_{n-1})}{(x-x'_1)\cdots(x-x'_n)}\\
\phantom{0=}{} +X_2(x)\sum_1^\infty c_n \eta_n\frac{(x-x_{-1})(x-x_0)\cdots (x-x_{n-2})}{(x-x'_1)\cdots(x-x'_{n-1})}\\
\phantom{0}{} =-c_0 c(x)-d(x)+c_1 X_2(x) \eta_1 (x-x_{-1})\\
\phantom{0=}{} +X_2(x)\sum_1^\infty (c_n \xi_n+c_{n+1}\eta_{n+1})\frac{(x-x_{-1})(x-x_0)\cdots (x-x_{n-2})(x-x_{n-1})}{(x-x'_1)\cdots(x-x'_n)}\\
\phantom{0}{} = (x-x_{-1})X_2(x)\left[ -c_0\beta-\delta+c_1\eta_1 +\sum_1^\infty (c_n \xi_n+c_{n+1}\eta_{n+1})\mathcal{X}_n(x)\right].
\end{gather*}
$X_2$ is a factor everywhere, from \eqref{factX2}, so
\begin{gather*}
0 = -c_0 (\beta x+\gamma) -(\delta x+\epsilon) +c_1 C_1 \displaystyle \frac{a(x_0)-c(x_0)(y_1-y_0)/2}{x_0-x'_1}
  (x-x_{-1})\\
\phantom{0 =}{}
+\sum_1^\infty (c_n \xi_n+c_{n+1}\eta_{n+1})\mathcal{X}_n(x),\\
 c_0=f(y_0)=  \frac{d(x_{-1})}{a(x_{-1})/(y_0-y_{-1})-c(x_{-1})/2}=-\frac{d(x_{-1})}{c(x_{-1})}
 =-(\delta x_{-1}+\epsilon)/(\beta x_{-1}+\gamma) ,\\
 c_1=  \frac{(\delta+\beta c_0)(x_0-x'_1)}{C_1(a(x_0)-c(x_0)(y_1-y_0)/2)}=
      \frac{(\gamma\delta-\beta\epsilon)(y_1-y'_1)X_2(x'_0)}{(y_1-y'_0)(x_0-x'_0)(a(x_0)-c(x_0)(y_1-y_0)/2)},
\end{gather*}
 as
 \begin{gather*}
 \frac{c_{n+1}}{c_n} = -\frac{\xi_n}{\eta_{n+1}} =
- \frac{C_n}{C_{n+1}} \frac{a(x'_n) +c(x'_n) (y'_{n+1}-y'_n)/2}{ a(x_n) -c(x_n) (y_{n+1}-y_n)/2}
\frac{(x_n-x_{-1})(x_n-x'_0)(x_n-x'_{n+1}) } { (x'_n-x_{-1})(x'_n-x'_0)(x'_n-x_{n-1}) },
\\
 c_n = \cdots \frac{ x'_n-x_{n-1}) }{C_n} \prod^{n-1}
     \frac{a(x'_k) +c(x'_k) (y'_{k+1}-y'_k)/2}{ a(x_k) +c(x_k) (y_{k+1}-y_k)/2}
  \mathcal{X}_n(x_{-1}) \mathcal{X}_n(x'_0).\tag*{\qed}
\end{gather*}\renewcommand{\qed}{}
\end{proof}

The formula \eqref{cn} achieves a construction of hypergeometric
type, as each term is a product of values of elliptic functions
with arguments in arithmetic progression. The exact order of each
term, i.e., the number of zeros and poles in a minimal
parallelogram, is not obvious \cite{SpZR}. Of course, a factor
like, say, $x_{-1}-x_k$ is an elliptic function of order 2 of
$t_0+kh$ from \eqref{ell12}. The same order holds for the ratio
\[
\frac{x_{-1}-x_k}{y_{-1}-y_k}=\frac{x_{-1}-\mathcal{E}_1(t_0+kh)}
{y_{-1}-\mathcal{E}_2(t_0+kh)},
\] as zeros of the numerator and the
denominator cancel each other.

Similar ef\/fects probably hold in other ratios encountered in \eqref{cn}, such as
\[
 \frac{a(x_k) -c(x_k) (y_{k+1}-y_k)/2}{(x_k-x_{-1})(x_k-x'_0)}
 \] but it is
not clear if more can be obtained by keeping elementary means, or if more elliptic
function machinery (theta functions) is needed.
An elementary description holds however in the ``logarithmic''
case $c(x)\equiv 0$. Then, \eqref{xp0}~and~\eqref{xm1} already
tell that $x_{-1}$ and $x'_0$ are two roots of $a(x)=0$. And as
the polynomial $a$ has degree~3 in Theorem~\ref{theorem1}, let
$a(x)=(x-x_{-1})(x-x'_0)(x-\zeta)$. Then, from \eqref{cn},
\begin{gather}
c_n = -c_1 \frac{C_1}{ x'_1\!-x_0 }  (x'_n\!-x_{n-1})
     \frac{  (y_{-1}\!-y'_1)\cdots (y_{-1}\!-y'_{n-1}) X_2(x_{-1}) (x_{-1}\!-x_0)\cdots(x_{-1}\!-x_{n-2})  }
           {  (y_{-1}\!-y_1) \cdots (y_{-1}\!-y_{n-2})  (x_{-1}\!-x'_0)\cdots (x_{-1}\!-x'_n)  } \nonumber\\
\phantom{c_n =}{}\times     \prod_{k=0}^{n-1}\frac{a(x'_k)} { a(x_k) }
              \frac{(x_k-x_{-1})(x_k-x'_0)}{(x'_k-x_{-1})(x'_k-x'_0)},\nonumber\\
  c_n \mathcal{Y}_n(x) = -c_1 \frac{C_1}{ x'_1-x_0 }  (x'_n-x_{n-1})\nonumber\\
 \phantom{c_n \mathcal{Y}_n(x)=}{}\times
     \frac{  (y_{-1}-y'_1)\cdots (y_{-1}-y'_{n-1}) X_2(x_{-1}) (x_{-1}-x_0)\cdots(x_{-1}-x_{n-2})  }
           {  (y_{-1}-y_1) \cdots (y_{-1}-y_{n-2})  (x_{-1}-x'_0)\cdots (x_{-1}-x'_n)  } \nonumber\\
\phantom{c_n \mathcal{Y}_n(x)=}{}\times      \prod_{k=0}^{n-1}\frac{(x'_k-\zeta)(x-y_k)}
                                                     { (x_k-\zeta)(x-y'_{k+1}) }.\label{cnlog}
\end{gather}

\section{A word on convergence}\label{conv}

\subsection{Average behaviour}

We expect products occurring in \eqref{cn} or \eqref{cnlog} to
behave like powers, like
\[
 \prod_1^n (x-x_k) =\prod_1^n (x-\mathcal{E}(t_0+kh))\approx \Phi_+(x)^n.
 \]
 What is $\Phi_+(x)=\exp\mathcal{V}_+(x)$, where $\mathcal{V}_+$ is the complex potential of the distributions of  $x_k$?
  For~$x'_k$, we
write $\mathcal{V}_-(x)$.  For~$y$, let us use the symbol
$\mathcal{W}$.

The main behaviour of the $n^{\text{th}}$ term of \eqref{cnlog} is
therefore
\begin{gather}
\exp\big(n(\mathcal{W}_-(y_{-1})-\mathcal{W}_+(y_{-1})+\mathcal{V}_+(x_{-1})
        -\mathcal{V}_-(x_{-1})+\mathcal{V}_-(\zeta)\nonumber\\
\qquad{}-\mathcal{V}_+(\zeta)
        +\mathcal{W}_+(x)-\mathcal{W}_-(x)  )\big).\label{pot1}
\end{gather}

Remark that we will only need $\mathcal{V}=\mathcal{V}_+ -\mathcal{V}_-$ and
$\mathcal{W}=\mathcal{W}_+ -\mathcal{W}_-$.

If $h$ is a general complex number, $x_k$ f\/ill the whole
complex plane and no convergence occurs.

Let $h$ be a  \emph{real}  irrational multiple of a period $\omega$,
then the same factors reappear approximately in the product after $N$ steps
if $Nh$ is close to an integer times $\omega$. $\Phi(x)$ is the limit
of the $N^{\text{th}}$ roots of such products. The various $kh$, for
$k=1,2,\dots, N$, modulo $\omega$, f\/ill uniformly the segment $[0,\omega]$, and
 $x_k$ f\/ill a curve which is the set of $\mathcal{E}(t_0+u), u\in[0,\omega]$:
for any $j$ in $\{1,2,\dots, N\}$, there is a $k$ such that $kh$ is close to
$j\omega/N$ modulo $\omega$. Indeed, let $Nh$ be close to $M_N\omega$, with gcd$(N,M_N)=1$. Then,
\[
 kh -\frac{j\omega}{N} = \omega\left(\frac{h}{\omega}-\frac{M_N}{N}\right)k
+\omega\frac{k M_N -j}{N},\] to any $j$, there are  integers $k$
and $m$ such that $k M_N -mN=j$ (Bezout).

 So, we rearrange the
product as
\[
\Phi(x) \sim
 \left[ \prod_{j=1}^N  (x-\mathcal{E}(j\omega/N +t_0))\right]^{1/N}
   \sim \exp\left[\frac{1}{\omega}\int_0^\omega \log(x-\mathcal{E}(u +t_0))\,du\right].
   \]
As $\mathcal{E}$ is the inversion of an elliptic integral of the f\/irst kind,
\[
 u+t_0= \int^{\mathcal{E}} \frac{dv}{\sqrt{P(v)}},
\]
we have
\[
\Phi(x)= \exp\left[\frac{1}{\omega}\int_{\{x_n\}}
   \frac{\log(x-v)\,dv}{\sqrt{P(v)}}\right],
   \]
where $\{x_n\}$ is the locus $=\{\mathcal{E}(u+t_0)\}$,
$u\in [0,\omega]$. The constant $1/\omega$ is such that
$\Phi(x)\sim x$ for large $x$:
\[
\omega=\int_{\{x_n\}}
   \frac{ dv}{\sqrt{P(v)}}.
   \]

 So, let the complex potential {\samepage
 \[
 \mathcal{V_+}(x) =
  \frac{1}{\omega} \int_{\{x_n\}} \frac{\log(x-v)\,dv}{\sqrt{P(v)}}
 \]
($\mathcal{V}_-$ will be used with  $x'_n$, and
$\mathcal{W}_{\pm}$ for  $y_n$ and $y'_n$).}

The formula for the potential will be linear after a convenient
conformal map.

One has the derivative
\[
 \mathcal{V'_+}(x) =
  \frac{1}{\omega} \int_{\{x_n\}} \frac{dv}{(x-v)\sqrt{P(v)}},
\]
with $\xi$ such that $x=\mathcal{E}(\xi)$, $dx/d\xi=\sqrt{P(x)}$.

 So, $\mathcal{V}'_+(x)$ and $\mathcal{V}'_-(x)$
are contour integrals on the locii f\/illed by $\{x_n\}$ and $\{x'_n\}$ drawn by
$\mathcal{E}(nh+t_0)$ and $\mathcal{E}(nh+t'_0)$. If $x$ is between these two locii,
the two contour integrals of $ \frac{dv}{(x-v)\sqrt{P(v)}}$ are the same
for $\mathcal{V}'_+(x)$ and $\mathcal{V}'_-(x)$, up to the residue at $v=x$:
\[ \mathcal{V}'(x)= \mathcal{V}'_+(x) -\mathcal{V}'_-(x)=
\frac{2\pi i }{\omega \sqrt{P(x)}} \ \Rightarrow \
\frac{d\mathcal{V}(x)}{d\xi} = \frac{2\pi i}{\omega}. \]

We see that the real part of $\mathcal{V}$ remains constant on
lines in the $\xi$-plane such that $d\xi/\omega$ is real, i.e.,
on parallel lines sharing the $\omega$-direction.

Remember that the step $h$ has been supposed to be a real multiple
of $\omega$, so the arguments in arithmetic progression of step
$h$ in the $\xi$-plane of the elliptic functions def\/ining a
sequence $x_n$, or $y_n$, etc. happen to draw parallel lines with
the $\omega$-direction! The real part of
$\mathcal{V}(\zeta)-\mathcal{V}(x_{-1})$ occurring in \eqref{pot1}
is therefore $2\pi/|\omega|$ times the distance between, say,
$\xi_\zeta$, if $\zeta$ is the value of the elliptic function at
$\xi_\zeta$, and the line leading to the $\{x_n\}$ sequence.

The remaining terms of \eqref{pot1} lead to a convergence
behaviour dominated by
\begin{gather} \label{ratec}
   \exp(-n\; \text{Im}\, 2\pi (\xi_x-\xi_\zeta)/\omega),
\end{gather}    where $\xi_x$ is sent to $x$ by the elliptic function.

Of course, convergence holds while $x$ is between the locus of
$x_n$ and the corresponding locus (equipotential line) containing
$\zeta$.

In a Jacobian setting, evaluation of \eqref{ratec} typically
involves $\exp(-n\pi K'/K)$, well known in Zolotarev problems
solutions and generalizations \cite{Gon69}.

Rate of approximation has already been related to potential problems by
Walsh \cite[chapter~9]{wal}, in papers and books going back to the 1930s! See also
Ganelius \cite{Ga}. For more
recent surveys and papers, the works of  Bagby~\cite{Bagby}, and by  Gon\v car and colleagues are
recommended \cite{Gon69,Gon2003,GR,Gon1992,Gon2006}.

It is quite remarkable that conf\/igurations of particles in statistical physics \cite{LouSpi1,LouSpi2,LouSpi3}
are described in the same way than zeros and poles of rational approximations
 \cite{Bagby,Gon69,Gon2003,GR,Gon1992,Gon2006,Meinguet2000,SaffTot,Sta}.

\subsection{Exceptional cases}

The properties of the irrational number relating the step $h$ to a
period $\omega$ must also be conside\-red~\cite{Spirconv}. Indeed,
\eqref{cnlog} contains a division by a factor $(y_{-1}-y_{n-2})$
which is the dif\/ference of the values of a function of period
$\omega$ at arguments dif\/fering by an integer multiple of $h$, so
that the result will be small whenever $(n-1)h$ is close to an
integer multiple of $\omega$. The dif\/ference will never vanish, as
$h/\omega$ is irrational, but could become VERY small inf\/initely
often. The set of irrational $h/\omega$ that could destroy the
convergence estimate above is fortunately of vanishing measure in
the set of real numbers, as shown by Hardy and Littlewood in~\cite{HardyL} (and reproduced by Lubinsky in \cite[pp.~854--855 and 871]{Lub04}).


\subsection*{Acknowledgments}

Many thanks to the organizers of the
workshop ``Elliptic Integrable Systems, Isomonodromy Problems, and Hypergeometric Functions'' (Hausdorf\/f Center for Mathematics, Bonn, July 2008),
to A.~Aptekarev, B.~Beckermann, A.C.~Matos, F.~Wielonsky,
of the Laboratoire Paul Painlev\'e UMR 8524,
Universit\'e de Lille 1, France, who organized their
$3^{\text{\`emes}}$  Journ\'ees Approxi\-ma\-tion on May 15--16, 2008.
Many thanks too to R.~Askey, L.~Haine, M.~Ismail, F.~Nijhof\/f, A.~Ronveaux, and, of course,
V.~Spiridonov and A.~Zhedanov for their preprints, interest, remarks, and kind words.
Many thanks to the referees for expert and careful reading, and kind words too.

This paper presents research results of the Belgian Programme on
Interuniversity Attraction Poles, initiated by the Belgian Federal
Science Policy Of\/f\/ice.

\pdfbookmark[1]{References}{ref}
\LastPageEnding

\end{document}